\documentclass[final,leqno]{siamltex}
\usepackage{amsfonts}
\usepackage{bm}
\usepackage{graphicx}

\newtheorem{remark}[theorem]{Remark}

\title{A multi-fidelity stochastic collocation method using locally improved reduced-order models}

\author{M. Raissi\thanks{Department of Mathematical Sciences, 
		George Mason University, Fairfax, VA 22030 email: {\tt mraissi@gmu.edu}.}
        \and P. Seshaiyer\thanks{Department of Mathematical Sciences, 
		George Mason University, Fairfax, VA 22030 email: {\tt pseshaiy@gmu.edu}.}}

\begin{document}

\maketitle

\begin{abstract}
Over the last few years there have been dramatic advances in our understanding of mathematical and computational models of complex systems in the presence of uncertainty. This has led to a growth in the area of uncertainty quantification as well as the need to develop efficient, scalable, stable and convergent computational methods for solving differential equations with random inputs. Stochastic Galerkin methods based on polynomial chaos expansions have shown superiority to other non-sampling and many sampling techniques. However, for complicated governing equations numerical implementations of stochastic Galerkin methods can become non-trivial. On the other hand, Monte Carlo and other traditional sampling methods, are straightforward to implement. However, they do not offer as fast convergence rates as stochastic Galerkin. Other numerical approaches are the stochastic collocation (SC) methods, which inherit both, the ease of implementation of Monte Carlo and the robustness of stochastic Galerkin to a great deal. However, stochastic collocation and its powerful extensions, e.g. sparse grid stochastic collocation, can simply fail to handle more levels of complication. The seemingly innocent Burgers equation driven by Brownian motion is such an example. In this work we propose a novel enhancement to stochastic collocation methods using locally improved deterministic model reduction techniques that can handle this pathological example and hopefully other more complicated equations like Stochastic Navier-Stokes. Local improvements to reduced-order models are achieved using sensitivity analysis of the proper orthogonal decomposition. Our numerical results show that the proposed technique is not only reliable and robust but also very efficient.
\end{abstract}



\pagestyle{myheadings}
\thispagestyle{plain}
\markboth{M. RAISSI AND P. SESHAIYER}{MULTI-FIDELITY STOCHASTIC COLLOCATION}

\section{Introduction}
The effectiveness of stochastic partial differential equations (SPDEs) in modelling complicated phenomena is a well-known fact. One can name wave propagation \cite{papanicolaou1971wave}, diffusion through heterogeneous random media \cite{papanicolaou1995diffusion}, randomly forced Burgers and NavierStokes equations (see e.g \cite{bensoussan1973equations,da2003ergodicity,khanin1997probability,mikulevicius2004stochastic} and the references therein) as a couple of examples. Currently, Monte Carlo is by far the most widely used tool in simulating models driven by SPDEs. However, Monte Carlo simulations are generally very expensive. To meet this concern, methods based on the Fourier analysis with respect to the Gaussian (rather than Lebesgue) measure, have been investigated in recent decades. More specifically, Cameron--Martin version of the Wiener Chaos expansion (see, e.g. \cite{cameron1947orthogonal,hida1993white} and the references therein) is among the earlier efforts. Sometimes, the Wiener Chaos expansion (WCE for short) is also referred to as the Hermite polynomial chaos expansion. The term polynomial chaos was coined by Nobert Wiener \cite{wiener1938homogeneous}. In Wieners work, Hermite polynomials served as an orthogonal basis. The validity of the approach was then proved in \cite{cameron1947orthogonal}. There is a long history of using WCE as well as other polynomial chaos expansions in problems in physics and engineering. See, e.g. \cite{crow1970relationship,orszag1967dynamical,chorin1971hermite,chorin1974gaussian}, etc. Applications of the polynomial chaos to stochastic PDEs considered in the literature typically deal with stochastic input generated by a finite number of random variables (see, e.g. \cite{sakamoto2002simulation,ghanem2003stochastic,xiu2003modeling,zhang2004efficient}). This assumption is usually introduced either directly or via a representation of the stochastic input by a truncated Karhunen--Lo\`{e}ve expansion. Stochastic finite element methods based on the Karhunen--Lo\`{e}ve expansion and Hermite polynomial chaos expansion \cite{ghanem2003stochastic,sakamoto2002simulation} have been developed by Ghanem and other authors. Karniadakis et al. generalized this idea to other types of randomness and polynomials \cite{jardak2002spectral,xiu2003modeling,xiu2002wiener}. The stochastic finite element procedure often results in a set of coupled deterministic equations which requires additional effort to be solved. To resolve this issue, stochastic collocation (SC) method was introduced. In this method one repeatedly executes an established deterministic code on a prescribed node in the random space defined by the random inputs. The idea can be found in early works such as \cite{mathelin2003stochastic,tatang1997efficient}. In these works mostly tensor products of one-dimensional nodes (e.g., Gauss quadrature) are employed. Tensor product construction despite making mathematical analysis more accessible (cf. \cite{babuvska2007stochastic}) leads to the curse of dimensionality since the total number of nodes grows exponentially fast as the number of random parameters increases. In recent years we are experiencing a surge of interest in the high-order stochastic collocation approach following \cite{xiu2005high}. The use of sparse grids from multivariate interpolation analysis, is a distinct feature of the work in \cite{xiu2005high}. A sparse grid, being a subset of the full tensor grid, can retain many of the accuracy properties of the tensor grid. While keeping high-order accuracy, it can significantly reduce the number of nodes in higher random dimensions. Further reduction in the number of nodes was pursued in \cite{agarwal2009domain,ma2009adaptive,nobile2008sparse,nobile2008anisotropic}. Applications of these numerical methods take a wide range. Here we mention some of the more representative works. It includes Burgers’ equation \cite{hou2006wiener,xiu2004supersensitivity}, fluid dynamics \cite{knio2006uncertainty,knio2001stochastic,le2002stochastic,lin2007stochastic,xiu2003modeling}, flow-structure interactions \cite{xiu2002stochastic}, hyperbolic problems \cite{chen2005uncertainty,gottlieb2008galerkin,lin2006predicting}, model construction and reduction \cite{doostan2007stochastic,ghanem2005identification,ghanem2006construction}, random domains with rough boundaries \cite{canuto2007fictitious,lin2007random,tartakovsky2006stochastic,xiu2006numerical}, etc.

Along with an attempt to reduce the number of nodes used by sparse grid stochastic collocation, one can try to employ more efficient deterministic algorithms. The current trend is to repeatedly execute a full-scale underlying deterministic simulation on prescribed nodes in the random space. However, model reduction techniques can be employed to create a computationally cheap deterministic algorithm that can be used for most of the grid points. This way we can limit the employment of an established while computationally expensive algorithm to only a relatively small number of points. A similar method is being used by K. Willcox and her team but in the context of optimization \cite{robinson2008surrogate}. ``Multifidelity", which we also adopt, is the term they employed in their work. For a tractable demonstration of our method we choose stochastic Burgers equation studied in \cite{hou2006wiener}. We believe that the non-linear nature of the Burgers equation with the extra complexity of Brownian motion, makes this equation a proper test for our method. The regular repetitive execution of the full-scale underlying deterministic algorithm for the sparse grid stochastic collocation simply fails for this equation. This equation is also studied in \cite{hou2006wiener} using Wiener Chaos expansion. For the deterministic Burgers equation resulting from stochastic collocation samples, we are going to use a high fidelity (accurate but computationally expensive) algorithm called ``group finite element'' (GFE) and a low fidelity (less accurate but computationally cheap) algorithm called ``group proper orthogonal decomposition''. These two methods are discussed in \cite{dickinson2010nonlinear}. The GFE method, also known as product approximation, is a finite element (FE) technique for certain types of non-linear partial differential equations. It expresses the non-linear terms of a PDE in a grouped form. As a result, spatial discretization of non-linear terms is computed only once before integration. Therefore, a substantial reduction in computational cost  is achieved \cite{christie1981product,fletcher1983group,roache1998verification}. Experiments with the GFE method have indicated an increase in economy and a slight increase in the nodal accuracy compared to FE solutions of the unsteady Burgers’ equations and many other problems \cite{christie1981product,fletcher1983group,srinivas1992computational}. Although theoretical results exist for other problems \cite{christie1981product,sanz1984interpolation,chen1989error,douglas1975effect,larsson1989interpolation,murdoch1986error,tourigny1990product}, the authors are unaware of convergence theory for the GFE applied to Burgers’ equation. The computational advantage of the GFE method over the conventional FE method for two and three dimensional Burgers’ equations and viscous compressible flows is demonstrated in \cite{fletcher1983group,srinivas1992computational}. As for the Group POD method, we use the projection of grouped non-linear terms onto a set of highly relevant global basis functions. This projection onto global basis functions further reduces the cost of simulation due to symmetry in the non-linear terms. Reduced order modelling, using proper orthogonal decompositions (POD) along with Galerkin projection, for fluid flows has seen extensive applications studied in \cite{sirovich1987turbulence,chambers1988karhunen,holmes1998turbulence,fahl2001computation,iollo2000stability,
kunisch2001galerkin,kunisch2002galerkin,henri2002stability,rowley2004model,camphouse2005boundary}. Proper orthogonal decomposition (POD) was introduce in Pearson \cite{pearson1901liii} and Hotelling \cite{hotelling1933analysis}. Since the work of Pearson and Hotelling, many have studied or used POD in a range of fields such as oceanography \cite{bjornsson1997manual}, fluid mechanics \cite{sirovich1987turbulence,holmes1998turbulence}, system feedback
control \cite{ravindran2000reduced,atwell2001reduced,atwell2001proper,kepler2002reduced,atwell2004reduced,lee2005reduced}, and system modeling \cite{fahl2001computation,kunisch2002galerkin,rowley2004model,henri2005convergence}.

The idea of our multi-fidelity stochastic collocation method is very simple. For each point in the stochastic parameter domain we search to see if the resulting deterministic problem is already solved by the high fidelity (GFE) algorithm for a sufficiently close problem. If yes, we use the solution to the nearby problem to create POD basis functions and we employ POD-Galerkin method (the low fidelity algorithm) to solve the original deterministic problem. However, a serious limitation of the POD method is that it is empirical. In other words, this basis accurately represents the data used to generate it, but may not be accurate when applied off-design. Thus, the reduced-order model may lose accuracy when applied to a problem with input data different from those used to generate the POD basis. To deal with this issue, we use sensitivity analysis in the basis selection step (see e.g. \cite{hay2009local,borggaard1997pde}) to locally improve the performance of POD basis functions. The newly derived basis allows for a more accurate solution of the problem when exploring the  stochastic parameter space.

\section{Stochastic Burgers - A tractable while serious test}
Our objective in this paper is to demonstrate how sparse grid stochastic collocation method can be enhanced using model reduction techniques. To simplify the presentation while employing a serious test, we study the following stochastic Burgers equation,
\begin{equation}
u_{t} + \frac{1}{2}(u^{2})_{x} = \mu u_{xx} + \sigma(x)\dot{W}(t), \label{burg}
\end{equation}
$(t,x) \in (0,T]\times[0,1], u(0,x) = u_{0}(x), u(t,0)=u(t,1)=0,$ where $W(t)$ is a Brownian motion and $u_{0} \in L_{2}([0,1])$ is a deterministic initial condition. If $\Vert u_{0}\Vert_{L_{2}} < \infty$ and $\Vert \sigma\Vert_{L_{2}} < \infty$, it is known (see, e.g \cite{da1994stochastic}) that (\ref{burg}) has a unique square integrable solution. If $u$ is a solution of equation (\ref{burg}), then $u$ is not only a function of $t$ and $x$, but it is also a function of the Brownian motion path $W^{t}_{0} = \{W(s), 0 \leq s \leq t\}$. In order to approximate the Brownian motion with a finite number of independent and identically distributed (iid) standard normal random variables, we let $\{h_{k}, k = 1,2,\ldots\}$ be an arbitrary orthonormal basis in $ L_{2}([0,t])$ and define $\xi_{k} := \int_{0}^{t}h_{k}(s)dW(s),$ for $k \in \{1,2,\ldots\}$. It can be shown that $\xi_{k}$ are iid Gaussian random variables. It is a standard fact that we can expand $W(s)$ as
\begin{equation}
W(s) = \int_{0}^{t}\chi_{[0,s]}(\tau)dW(\tau) = \sum_{k=1}^{\infty}\xi_{k}\int_{0}^{s}h_{k}(\tau)d\tau, \label{expansion}
\end{equation}
where $\chi_{[0,s]}(\tau)$ is the characteristic function of the interval $[0,s]$. Note that
\[
\chi_{[0,s]}(\tau) = \sum_{k=1}^\infty c_k h_k(\tau),
\]
where $c_k = \int_0^t \chi_{[0,s]}(\tau)h_k(\tau)d\tau = \int_0^s h_k(\tau)d\tau$. If $\{h_{k}, k = 1,2,\ldots\}$ are chosen as Haar wavelets, then expansion (\ref{expansion}) is exactly the Levy-Ciesielski construction \cite{mckean1969stochastic} of Brownian motion. Therefore we can view a solution $u$ of (\ref{burg}) as a function of $x, t$ and $\bm\xi = (\xi_{1},\xi_{2}, \ldots)$. We also know that expansion (\ref{expansion}) converges in the mean square sense (i.e., $E[W(s) - \sum_{k=1}^d\xi_{k}\int_0^s h_{k}(\tau)d\tau]^2$ goes to zero as  $d\rightarrow\infty$ uniformly for $s \leq t$.). Our first step is to truncate expansion (\ref{expansion}) at some point $d$ to get
\begin{equation}
u_{t} + \frac{1}{2}(u^{2})_{x} = \mu u_{xx} + \sigma(x)\sum_{k=1}^{d}\xi_{k}h_{k}(t).\label{truncated}
\end{equation}
In the followings, we are seeking an approximation to a random field $u_{d}(\bm{\hat{\xi}};t,x) \in C^{\infty}(\mathbb{R}^{d};L_{2}([0,T];W([0,1])))$ which satisfies equation (\ref{truncated}) along with the initial and boundary conditions of problem (\ref{burg}), where $\bm{\hat{\xi}} = (\xi_{1},\xi_2, \ldots,\xi_d)$ and $W([0,1])$ is a Banach space of functions $v:[0,1]\longrightarrow\mathbb{R}$. Let us assume that the random vector $\bm{\hat{\xi}}$ has $\rho(\bm{\hat{\xi}})$ (which is simply the multiplication of probability densities of $d$ standard normal variables) as its joint probability density. We are also interested in approximating expectation and higher moments $E[u_d^{k}] \in L_{2}([0,T];W([0,1])), k = 1,2,\ldots$ of the solution which are given by
\begin{equation}
E[u_d^k](t,x) = \int_{\mathbb{R}^d}u_d^{k}(\bm{\hat{\xi}};t,x)\rho(\bm{\hat{\xi}})d\bm{\hat{\xi}},~~k \in \{1,2,\ldots\}.\label{Moments}
\end{equation}
For simplicity, we are misusing the symbol $\bm{\hat{\xi}}$ to denote both a random vector and its realizations. For each fixed $\bm{\hat{\xi}}$, equation (\ref{truncated}) can be written as,
\begin{equation}
u_{t} + \frac{1}{2}(u^{2})_{x} = \mu u_{xx} + f_d(t,x),\label{burgers}
\end{equation}
where $f_d(t,x) = \sigma(x)\sum_{k=1}^{d}\xi_{k}h_{k}(t)$. Whether we are using Monte Carlo or Stochastic Collocation, in order to solve the stochastic equation (\ref{truncated}), we need to solve the deterministic equation (\ref{burgers}) many times for different values of $\bm{\hat{\xi}}$. In the followings, we specify the high fidelity algorithm for equation (\ref{burgers}) and the low fidelity algorithm for
\begin{equation}
u_{t} + \frac{1}{2}(u^{2})_{x} = \mu u_{xx} + g_d(t,x),\label{burgersNew}
\end{equation}
where $g_d(t,x) = \sigma(x)\sum_{k=1}^{d}\zeta_{k}h_{k}(t)$, which is practically equation (\ref{burgers}) evaluated at the new parameter $\bm{\hat{\zeta}} = (\zeta_1,\ldots,\zeta_d)$. Note that $\bm{\hat{\xi}}$ and $\bm{\hat{\zeta}}$ are not necessarily equal. Consequently, $f_d$ and $g_d$ are generally different.

\section{GFE as a high-fidelity deterministic algorithm}
A standard finite element approximation to the solution $u_d(t,x)$ of equation (\ref{burgers}) can be written as $u_{d,N}(t,x) = \sum_{j=1}^{N}\alpha_{j}(t)\beta_j(x)$, where $\{\beta_j(x), j = 1,\ldots,N\}$ are $N$ piecewise linear finite element basis functions, and each $\alpha_j(t)$ is an unknown function of time. We also define $W_N([0,1]):=\textrm{span}\{\beta_j,j=1,\ldots,N\}$. The basis functions are chosen according to a computational grid on the domain [0,1] with $\{x_n, n= 0,\ldots,N+1\}$ as its set of nodes, in a way that $\beta_j(x_n) = \delta_{j,n}$, where $\delta_{j,n}$ is the Kronecker delta and $j,n = 1,\ldots,N$. Note that $u^2_{d,N}(t,x_n) = \sum_{j=1}^N\alpha_j^2(t)\beta_j(x_n)$. This motivates us to approximate $u_d^2(t,x)$ by $\sum_{j=1}^N\alpha_j^2(t)\beta_j(x)$. Use of the weak form of (\ref{burgers}) and the approximations for $u_d$ and $u_d^2$ results in the following differential equations.
\begin{eqnarray}
\bm M\dot{\bm{\alpha}} &=& -\mu\bm A\bm\alpha -\frac{1}{2} \bm G(\bm\alpha)+\bm V(t),\label{GFE}\\
\bm \alpha(0) &=& \bm \alpha_{\bm 0} = \bm M^{-1}[(u_0,\beta_i)]_{i=1}^{N},\nonumber
\end{eqnarray} 
where $\bm G(\bm \alpha) = \bm N[\textrm{diag}(\bm \alpha)]\bm\alpha$, $[\bm N]_{ij}=(\beta'_j,\beta_i)$, $[\bm M]_{ij} = (\beta_j,\beta_i)$, $[\bm A]_{ij} = (\beta'_j,\beta'_i)$, $[\bm V(t)]_i = (f_d(t,.),\beta_i)$ and $(f,g) = \int_0^1f(x)g(x)dx$ is the standard $L_2([0,1])$ inner product. We are using $\dot{y}$ and $y'$ to denote the derivatives with respect to time and space, respectively. This algorithm is called group finite elements (GFE). A detailed study of this method can be found in \cite{dickinson2010nonlinear}.

\section{Group POD as a low-fidelity deterministic algorithm} Solution to equations (\ref{GFE}), with $\bm V(t) = [(f_d(t,.),\beta_i)]_{i=1}^N$, gives an approximate solution to (\ref{burgers}) given by $w(t,x)=\sum_{j=1}^N\alpha_j(t)\beta_j(x)$. Let $\{w(t_i,.),i=1,\ldots,S\}$ be a set of $S$ ``snapshots'', where $t_1<t_2<\ldots<t_S$ are equally spaced points of time in the interval $[0,T]$. Furthermore, define $\bm{Y} = [\bm\alpha(t_1) \cdots \bm\alpha(t_S)]$ to be the snapshots matrix. Let $\bm L$ be a lower triangular matrix resulting from the Cholesky decomposition of the mass matrix $\bm M$, i.e. $\bm M = \bm L \bm L^T$. Define $\bm{\tilde{Y}}:= \bm L ^T \bm Y$ to be the weighted snapshots matrix. The correlation matrix $\bm K$ of the data set $\{w(t_i,.),i=1,\ldots,S\}$ can be defined as
\begin{equation}
\bm K := \left(\frac{1}{S}(w_i,w_j)\right)_{i,j=1}^S = \frac{1}{S}\bm Y^T\bm M \bm Y = \frac{1}{S}\bm{\tilde{Y}}^T \bm{\tilde{Y}}\label{CorrelationMatrix}
\end{equation}
where $w_i = w(t_i,.),i=1,\ldots,S$. Let $\{\lambda_k,\mathcal{Z}_k\}$ denote the eigenvalues and the corresponding normalized eigenvectors of $\bm K$. Define $\bm{\mathcal{Z}}$ to be the matrix $[\mathcal{Z}_1| \ldots |\mathcal{Z}_S]$. The POD basis functions $\{\psi_k \}_{k=1}^S$ are given as
\begin{equation}
\psi_k = \frac{1}{\sqrt{S\lambda_k}}\sum_{i=1}^S[\bm{\mathcal{Z}}]_{i,k}w_i, ~~~k = 1,\ldots,S \label{PODbasis}.
\end{equation}
Note that
\begin{equation}
\frac{1}{S}\bm{\tilde{Y}}^T \bm{\tilde{Y}} \bm{\mathcal{Z}} = \bm{\mathcal{Z}} \bm{\Lambda},
\end{equation}
where $\bm{\Lambda} := diag(\lambda_k,k=1,\ldots,S)$. Let us express the POD basis function $\psi_j$ as $\psi_j(x) = \sum_{i=1}^N\psi_{i,j}\beta_i(x)$ and let $\bm \psi = (\psi_{i,j})$ for $i=1,\ldots,N$ and $j=1,\ldots,S$. Using (\ref{PODbasis}) we get
\begin{equation}
\bm \psi = \bm Y \bm{\mathcal{Z}} (S \bm \Lambda)^{-\frac{1}{2}}.\label{PODmodes}
\end{equation}

Let the solution $u_d$ to (\ref{burgersNew}) be written as
\begin{equation}
u_d(t,x) = U(x) + v_d(t,x)\label{Averaged}
\end{equation}
where $U(x):=\frac{1}{S}\sum_{i=1}^{S}w(t_i,x)$. Therefore, using (\ref{burgersNew}), we obtain that $v_d$ satisfies
\begin{equation}
v_t + \frac{1}{2}(v^2)_x -  \mu (v_{xx}+U'')+ UU' + Uv_x + vU' =g_d,\label{GPOD}
\end{equation}
$(t,x)\in (0,T]\times[0,1], v(0,x) = v_0(x) = u_0(x) - U(x), v(t,0)=v(t,1)=0$.
Let the POD approximation of $v_d(t,x)$ be given by $v_{d,p}(t,x) = \sum_{j=1}^Ma_j(t)\psi_j(x)$, where $\psi_j$ is a POD basis function, and $a_j(t)$ is an unknown function of time.  Furthermore, let $v_d^2(t,x)$ be approximated by $\sum_{j=1}^MF_j(\bm a)\psi_j(x)$, with $F_j(\bm a)$ as an unknown function of $\bm a(t) = [a_j(t)]_{j=1}^M$.
The projection of equation (\ref{GPOD}) onto the POD set ${\bm \Psi} = \{\psi_k \}_{k=1}^M$ results in the variational problem of finding $v \in L_2([0,T];W_{N,M}([0,1]))$, where $W_{N,M}([0,1]) := \textrm{span}\{\bm \Psi\}$, such that
\begin{eqnarray}
(v_t,\psi_i) + \frac{1}{2}((v^2)_x,\psi_i) + (UU',\psi_i) &+& \ldots\label{GPODweak}\\
\mu(v_x+U',\psi_i') &+& (Uv_x + U'v,\psi_i) = (g_d,\psi_i),\nonumber
\end{eqnarray}
and $(v(0,.),\psi_i) = (v_0,\psi_i)$, for $i = 1,\ldots,M$. When the POD approximations for $v_d$ and $v_d^2$ are substituted into equation (\ref{GPODweak}), we obtain the system of ordinary differential equations
\begin{eqnarray}
\bm M\dot{\bm{a}} &=& -\bm A\bm a - \frac{1}{2}\bm{N_p}\bm F(\bm a)-\bm V(t),\label{GPODode}\\
\bm a(0) &=& \bm a_{\bm 0} = [(v_0,\psi_i)]_{i=1}^{M},\nonumber
\end{eqnarray}
where $[\bm{N_p}]_{ij} = (\psi_j',\psi_i)$, $[\bm M]_{ij} = (\psi_j,\psi_i)$, $[\bm A]_{ij} = \mu(\psi'_j,\psi'_i)+(U\psi_j'+U'\psi_j,\psi_i)$, $[\bm V]_i = (-g_d+UU',\psi_i)+\mu(U',\psi_i')$, and $\bm F(\bm a) = [F_i]_{i=1}^{M}$ is determined as given bellow. At the grid points $x_n,~n=1,\ldots,N$ we want to have
\begin{equation}
\sum_{j=1}^MF_j(\bm a)\psi_j(x_n) = \left(\sum_{j=1}^Ma_j(t)\psi_j(x_n)\right)^2.
\end{equation}
Let $\gamma_{nj} = \psi_j(x_n)$, then,
\begin{equation}
\sum_{j=1}^MF_j(\bm a)\gamma_{nj} = \left(\sum_{j=1}^Ma_j(t)\gamma_{nj}\right)^2 = \sum_{j,\ell=1}^{M}\gamma_{nj}\gamma_{n\ell}a_ja_\ell.
\label{Nonlinear}
\end{equation}
Since $\gamma_{nj}\gamma_{n\ell}=\gamma_{n\ell}\gamma_{nj}$, for $j\neq\ell$ we can avoid computing $\gamma_{nj}\gamma_{n\ell}$ twice by writing
\begin{equation}
\sum_{j,\ell=1}^{M}\gamma_{nj}\gamma_{n\ell}a_ja_\ell = \hat{\bm\gamma}_n\hat{\bm a}
\end{equation}
where $\hat{\bm\gamma}_n$ is a $1\times\frac{1}{2}(M^2+M)$ vector given as
\begin{eqnarray}
\hat{\bm\gamma}_n := [\gamma_{n1}\gamma_{n1},~2\gamma_{n1}\gamma_{n2},\ldots,2\gamma_{n1}\gamma_{nM},~\gamma_{n2}\gamma_{n2},~2\gamma_{n2}\gamma_{n3},\ldots,2\gamma_{n2}\gamma_{nM},\ldots,\gamma_{nM}\gamma_{nM}]\nonumber
\end{eqnarray}
and the $\frac{1}{2}(M^2+M)\times 1$ vector $\hat{\bm a}$ is 
\begin{eqnarray}
\hat{\bm a} := [a_1a_1,a_1a_2,\ldots,a_1a_M,a_2a_2,a_2a_3,\ldots,a_2a_M,\ldots,a_Ma_M]^T\nonumber.
\end{eqnarray}
We may write (\ref{Nonlinear}) in the matrix form as
\begin{equation}
\bm\Gamma\bm F(\bm a) = \hat{\bm\Gamma}\hat{\bm a}\label{Nonlinear2}
\end{equation}
where the $N\times M$ matrix $\bm\Gamma$ is given by
\begin{equation}
\bm\Gamma = \left[ \begin{array}{cccc}
\gamma_{11} & \gamma_{12} & \ldots & \gamma_{1M} \\ 
\vdots & \vdots & \vdots & \vdots \\ 
\gamma_{N1} & \gamma_{N2} & \ldots & \gamma_{NM}
\end{array} \right]
\end{equation}
and $\hat{\bm\Gamma}$ is the $N\times \frac{1}{2}(M^2+M)$ matrix $[\hat{\bm\gamma}_i]_{i=1}^{N}$. Equation (\ref{Nonlinear2}) normally has no solution, since $M$ is generally chosen to be less than $N$. However, similar to the method used in \cite{dickinson2010nonlinear}, one can directly approximate the product $\bm{N_p}\bm F(\bm a)$ by $\bm\Gamma^T \bm N\bm\Gamma\bm F(\bm a)$, where $\bm N$ is the matrix used in (\ref{GFE}). We can finally approximate $\bm{N_p}\bm F(\bm a)$ by
\begin{equation}
\bm{N_p}\bm F(\bm a) \approx\hat{\bm N}\hat{\bm a},
\end{equation}
where $\hat{\bm N}$ is the $M\times\frac{1}{2}(M^2+M)$ matrix given by $\hat{\bm N} = {\bm \Gamma}^T\bm N\hat{\bm \Gamma}$. With this, problem (\ref{GPODode}) takes the form
\begin{eqnarray}
\bm M\dot{\bm{a}} &=& -\bm A\bm a - \hat{\bm{N}}\hat{\bm a}-\bm V(t),\label{GPODodeMain}\\
\bm a(0) &=& \bm a_{\bm 0} = [(v_0,\psi_i)]_{i=1}^{M}.\nonumber
\end{eqnarray}
This model is then solved for $\bm a(t)$ to give an approximation $v_{d,p}(t,x) = \sum_{j=1}^Ma_j(t)\psi_j(x)$ of $v_d(t,x)$. Therefore, according to (\ref{Averaged}), $u_d(t,x)$ is approximated by $U(x) + v_{d,p}(t,x)$.  The above algorithm called group proper orthogonal decomposition (group POD) is studied in details in \cite{dickinson2010nonlinear}.

\section{Local improvements to POD basis functions} Note that we want to use the POD basis functions found at $\bm{\hat{\xi}}$ to solve the deterministic problem (\ref{burgersNew}) at a nearby point $\bm{\hat{\zeta}}$. It is a well known fact that if $\bm{\hat{\zeta}} = \bm{\hat{\xi}}$ we get sufficiently accurate approximations of the solution of problem (\ref{burgersNew}) or equivalently (\ref{burgers}). However, for $\bm{\hat{\zeta}} \neq \bm{\hat{\xi}}$, we usually experience a drop in accuracy. In order to deal with issue, we use the sensitivity analysis of POD basis functions to improve the accuracy of the solution found at $\bm{\hat{\zeta}}$. Note that the POD bases (\ref{PODbasis}) are a function of $\bm{\hat{\xi}}$. We are specifically interested in the sensitivity of these POD bases in the direction of $\bm{\hat{\zeta}}-\bm{\hat{\xi}}$. To do so, we introduce a new parameter $\theta\in [0,1]$ and consider the mapping $\bm{\hat{\xi}}+\theta(\bm{\hat{\zeta}}-\bm{\hat{\xi}})$. Therefore, the POD bases can be considered to be a function of $\theta$. We start by noting that:
\begin{equation}
\bm K \mathcal{Z}_k = \lambda_k  \mathcal{Z}_k,\label{Eigenvalue}
\end{equation}
where $\bm K$ is defined by (\ref{CorrelationMatrix}), and $\{\lambda_k,\mathcal{Z}_k\}$ denote the eigenvalues and the corresponding normalized eigenvectors of $\bm K$. We assume that the entries of $\bm{K}$, $\bm{\mathcal{Z}}$, and $\bm{\Lambda}$ are smooth functions of the parameter $\theta$ so that (\ref{Eigenvalue}) can be differentiated with respect to $\theta$. In what follows, partial derivative of any matrix or vector is denoted using the superscript $(\theta)$. Therefore, by implicit differentiation of (\ref{Eigenvalue}) with respect to $\theta$ we get:
\begin{eqnarray}
(\bm K-\lambda_k \bm I)\mathcal{Z}_k^{(\theta)} &=& -(\bm K^{(\theta)}-\lambda_k^{(\theta)} \bm I)\mathcal{Z}_k.\label{eqn:PODSensitivity}
\end{eqnarray}
Equation (\ref{eqn:PODSensitivity}) has a solution only if the right-hand side vector belongs to the rang of $\bm{K}-\lambda_k \bm{I}$ and thus must be orthogonal to $ker(\bm{K}-\lambda_k \bm{I})$ which is spanned by $\mathcal{Z}_k$. Therefore, we should have that
\begin{equation}
\mathcal{Z}_k^T(\bm K^{(\theta)}-\lambda_k^{(\theta)} \bm I)\mathcal{Z}_k=0.
\end{equation}
Since $\mathcal{Z}_k$ has unit norm, the sensitivity of the eigenvalues is obtained by (see e.g. \cite{lancaster1964eigenvalues,fox1968rates,murthy1988derivatives,seyranian1994multiple})
\begin{equation}
\lambda_k^{(\theta)} = \mathcal{Z}_k^T\bm K^{(\theta)}\mathcal{Z}_k.
\end{equation}
Note that $\bm K = \frac{1}{S}\bm{\hat{Y}}^T\bm{\hat{Y}}$ and $\bm{\hat{Y}} = \bm{L}^T\bm Y$. Therefore, $\bm K^{(\theta)} = \frac{1}{S}(\bm{\hat{Y}}^{(\theta)^T}\bm{\hat{Y}} + \bm{\hat{Y}}^T\bm{\hat{Y}}^{(\theta)})$ and $\bm{\hat{Y}}^{(\theta)} = \bm L ^T \bm Y^{(\theta)}$. Provided that we know what $\bm Y^{(\theta)}$ should be, we can now fully characterize the solution $\mathcal{Z}_k^{(\theta)}$ of equation (\ref{eqn:PODSensitivity}). We find a particular solution $\mathcal{S}_k$ of (\ref{eqn:PODSensitivity}) in the least-square sense (obtaining the minimum norm solution). Since $\lambda_k$ is simple, for all $\varrho\in \mathbb{R}$, $\mathcal{S}_k + \varrho\mathcal{Z}_k$ is the general expression for the solutions of (\ref{eqn:PODSensitivity}). To determine the particular solution of (\ref{eqn:PODSensitivity}) which corresponds to the sensitivity $\mathcal{Z}_k^{(\theta)}$ of $\mathcal{Z}_k$, we need an additional condition. This comes naturally from the normalization condition $\mathcal{Z}_k^T\mathcal{Z}_k=1$ which was employed to specify $\mathcal{Z}_k$. Differentiating the normalization condition we get $\mathcal{Z}_k^T\mathcal{Z}_k^{(\theta)}=0$ and consequently $\varrho = - \mathcal{Z}_k^T\mathcal{S}_k$. Finally,
\begin{equation}
\mathcal{Z}_k^{(\theta)}=\mathcal{S}_k-(\mathcal{Z}_k^T\mathcal{S}_k)\mathcal{Z}_k.
\end{equation}
Once the sensitivity of matrices $\bm{\mathcal{Z}}$ and $\bm \Lambda$ are determined, the sensitivity of POD basis modes $\bm \psi$ and POD bases $\psi_k$ are straightforward to be computed by differentiating (\ref{PODmodes}). More specifically,
\begin{eqnarray}
\bm \psi ^{(\theta)} &=& \bm Y^{(\theta)}\bm{\mathcal{Z}}(S\bm{\Lambda})^{-1/2} + \bm Y\bm{\mathcal{Z}}^{(\theta)}(S\bm{\Lambda})^{-1/2}+\bm Y\bm{\mathcal{Z}}S^{-1/2}(\bm{\Lambda}^{-1/2})^{(\theta)}\nonumber\\
&=& \bm Y^{(\theta)}\bm{\mathcal{Z}}(S\bm{\Lambda})^{-1/2} + \bm Y\bm{\mathcal{Z}}^{(\theta)}(S\bm{\Lambda})^{-1/2}-\frac{1}{2}\bm Y\bm{\mathcal{Z}}S^{-1/2}(\bm{\Lambda}^{-1/2}\bm\Lambda^{(\theta)}\bm\Lambda^{-1})\nonumber\\
&=& \bm Y^{(\theta)}\bm{\mathcal{Z}}(S\bm{\Lambda})^{-1/2} + \bm Y\bm{\mathcal{Z}}^{(\theta)}(S\bm{\Lambda})^{-1/2}-\frac{1}{2}\bm\psi\bm\Lambda^{(\theta)}\bm\Lambda^{-1}.
\end{eqnarray}
To complete the sensitivity analysis of the POD bases, we still need to find the sensitivity of the snapshot data matrix $\bm Y^{(\theta)}$. This can be done by the sensitivity analysis of equation (\ref{burgers}).
\subsection{Sensitivity analysis of the Burgers equation} Let $u$ be the solution of equation (\ref{burgers}) and let $z = \partial_\theta u$ be the partial derivative of $u$ with respect the parameter $\theta$. Taking the derivative of equation (\ref{burgers}) with respect to $\theta$, we get that $z$ should satisfy
\begin{equation}
z_t + (zu)_x=\mu z_{xx}+\partial_\theta f_d(t,x),\label{eqn:Sensitivity}
\end{equation}
where $\partial_\theta f_d(t,x) = \partial_\theta[ \sigma(x)\sum_{k=1}^d(\xi_k+\theta(\zeta_k-\xi_k))h_k(t)] = \sigma(x)\sum_{k=1}^d(\zeta_k-\xi_k)h_k(t)$, with zero boundary and initial conditions. Note that equation (\ref{eqn:Sensitivity}) is no more non-linear and can be efficiently solved using a regular Finite Element method in a negligible amount of time. The finite element basis functions are assumed to be as before and we are seeking the solution in the space $W_N([0,1]) = \textrm{span}\{\beta_j, j=1,\ldots,N\}$ defined earlier. Let $z(t,x) = \sum_{i=1}^Nz_i(t)\beta_i(x)$ be the solution resulting from solving (\ref{eqn:Sensitivity}) with the Finite Element method. Then the sensitivity of the snapshot matrix $\bm Y$ is given by $\bm Y^{(\theta)} = (z_i(t_j))$, $i = 1,\ldots,N$ and $j=1,\ldots,S$. This completes the sensitivity analysis of POD basis functions in the directions of $\bm{\hat{\zeta}}-\bm{\hat{\xi}}$.
\subsection{Improving POD bases} Following \cite{hay2009local}, we state two ideas for constructing improved reduced bases.
\begin{enumerate}
\item \textit{Extrapolated basis:} Note that POD bases introduce in (\ref{PODbasis}) are functions of $\bm{\hat{\xi}}$ and consequently functions of $\theta$ used in the transformation $\bm{\hat{\xi}}+\theta(\bm{\hat{\zeta}} - \bm{\hat{\xi}})$. Let us use $\psi_k(\theta)$ to emphasize this dependence. Note that when $\theta = 0$ we are considering the POD basis functions at $\bm{\hat{\xi}}$ and when $\theta = 1$ we are considering them at $\bm{\hat{\zeta}}$. Now let us use
\begin{equation}
\psi_k(\theta) \simeq \psi_k(0) + \Delta \theta \frac{\partial \psi_k}{\partial \theta}(0)
\end{equation}
to approximate $\psi_k(\theta)$. The capability of this extrapolation obviously depends on the assumption that POD modes behave nearly linear with respect to the parameter $\theta$. However, using this method, the dimension of the reduced basis is preserved. Other approaches based on this idea can be found in \cite{lehmann2005wake,lieu2006reduced,morzynski2007continuous}. These papers motivate the extrapolation approach based on mode sensitivity, by showing an increase in robustness of the derived POD models with respect to parameter change.
\item \textit{Expanded basis:} The sensitivity of the modes generally seem to span a different subspace than the POD modes. Therefore, it seems plausible to expect that if we seek the approximate solution in the space spanned by the union of these two sets, we can represent a broader range of solutions. Therefore, we use $\{\partial_\theta\psi_k\}$ to expand the original POD basis functions $\{\psi_k\}$. By a misuse of notation, we are still using
\[W_{N,M}([0,1]) = \textrm{span}_{k=1,\ldots,M}\{\psi_k,\partial_\theta\psi_k\}\]
to denote the space where we seek the solution of equation (\ref{GPODweak}). The underlying assumption of this approach is that $W_{N,M}([0,1])$ is well suited to address the change in the solution induced by a change in parameter. This indeed is a legitimate assumption since the sensitivities represent changes in the parameter space. However,  the dimension of the reduced basis is now doubled.
\end{enumerate}
\section{Multi-fidelity stochastic collocation} We are finally in a position to demonstrate the multi-fidelity stochastic collocation method. Our aim to approximate the solution of (\ref{truncated}) in the space $V_{N,\mathbf{p}} = L_{2}([0,T];W_N([0,1]))\otimes\mathcal{P}_{\mathbf{p}}(\mathbb{R}^d)$, where $W_N([0,1]) = \textrm{span}\{\beta_j\}$ is the finite element space, and $\mathcal{P}_{\mathbf{p}}(\mathbb{R}^d)$ is the span of tensor product polynomials with degree at most $\mathbf{p} = (p_1,\ldots,p_d)$. Choose $\eta>0$ to be a small real number. The procedure for approximating the solution of (\ref{truncated}) is divided into two parts:
\begin{enumerate}
\item Fix $\bm{\hat{\zeta}}\in\mathbb{R}^d$, and search the $\eta$-neighbourhood $B_\eta(\bm{\hat{\zeta}})$ of $\bm{\hat{\zeta}}$. We use
$
B_\eta(\bm{\hat{\zeta}}) = \{\bm{\hat{\xi}}\in\mathbb{R}^d:|\xi_k-\zeta_k|<\eta, ~\forall k=1,\ldots,d\}.
$
If problem (\ref{burgers}) is not already solved by the GFE method for any nearby problem with $\bm{\hat{\xi}} \in B_\eta(\bm{\hat{\zeta}})$, let $\bm{\hat{\xi}} = \bm{\hat{\zeta}}$ and solve problem (\ref{burgers}) using the GFE method. In contrast, if equation (\ref{burgers}) is already solved for some points in $B_\eta(\bm{\hat{\zeta}})$, pick the closest one to $\bm{\hat{\zeta}}$ and call it $\bm{\hat{\xi}}$. In either case, use the solution at $\bm{\hat{\xi}} \in B_\eta(\bm{\hat{\zeta}})$ to find a small number of suitable basis functions resulting from local improvements to POD bases using sensitivity analysis. Let $W_{N,M}([0,1])\subset W_{N}([0,1])$ be the span of these basis functions. Now use the Group POD method to solve problem (\ref{burgersNew}) at $\bm{\hat{\zeta}}\in\mathbb{R}^d$ and get the solution $u_d(\bm{\hat{\zeta}},t,x)$.
\item Collocate on zeros of suitable orthogonal polynomials and build the interpolated solution $u_{d,\mathbf{p}}\in V_{N,\mathbf{p}}$ using
\begin{eqnarray}
& &u_{d,\mathbf{p}}(\bm{\hat{\zeta}},t,x) = \mathcal{I}_\mathbf{p}u_{d}(\bm{\hat{\zeta}},t,x) =\label{CollocationScheme}\\
& &\sum_{j_1=1}^{p_1+1}\cdots\sum_{j_d=1}^{p_d+1}u_d(\zeta_1,\ldots,\zeta_d,t,x)(l_{j_1}(\bm{\hat{\zeta}})\otimes\cdots\otimes l_{j_d}(\bm{\hat{\zeta}}),\nonumber
\end{eqnarray}
where the functions $\{l_{j_k}\}_{k=1}^d$ can be taken as Lagrange polynomials.  Using this formula, as described in \cite{babuvska2007stochastic}, mean value and variance of $u_d$ can also be easily approximated.
\end{enumerate}
\subsection{Generalization to sparse grid}
Here, we give a short description of the isotropic Smolyak algorithm. More detailed information can be found in \cite{barthelmann2000high,nobile2008sparse}. Assume $p_1=p_2=\cdots=p_d=p$. For $d=1$, let $\{\mathcal{I}_{1,i}\}_{i=1,2,\ldots}$ be a sequence of interpolation operators given by equation (\ref{CollocationScheme}). Define $\Delta_0 = \mathcal{I}_{1,0} = 0$ and $\Delta_i = \mathcal{I}_{1,i} - \mathcal{I}_{1,i-1}$. Now for $d > 1$, let
\begin{eqnarray}
\mathcal{A}(q,d) = \sum_{0\leq i_1+i_2+\ldots+i_d \leq q}\Delta_{i_1}\otimes\cdots\otimes \Delta_{i_d}
\end{eqnarray}
where $q$ is a non-negative integer. $\mathcal{A}(q,d)$ is the Smolyak operator, and $q$ is known as the sparse grid \textit{level}. Now instead of (\ref{CollocationScheme}), $\mathcal{A}(q,d)u_{d}(\bm{\hat{\zeta}},t,x)$ can be used to approximate the solution $u_d$ of (\ref{truncated}). This way one reduces the number of grid points on which the deterministic algorithms should be employed.
\section{Numerical Experiments}
In the followings, we consider equation (\ref{truncated}) given once again bellow for conveniences.
\begin{equation}
u_{t} + \frac{1}{2}(u^{2})_{x} = \mu u_{xx} + \sigma(x)\sum_{k=1}^{d}\xi_{k}h_{k}(t),\label{truncatedOnceAgain}
\end{equation}
$(t,x) \in (0,T]\times[0,1], u(0,x) = u_{0}(x), u(t,0)=u(t,1)=0,$ where $u_{0} \in L_{2}([0,1])$ is a deterministic initial condition. We assume that $\sigma(x) = 0.1\cos(4\pi x)$ and $\mu = 1/100$. We let $T = 0.8$ and $u_0(x) = (e^{\cos(5\pi x)} - \frac{3}{2})\sin(\pi x)$. We project the Brownian motion in $[0,T]$ on the trigonometric basis functions $h_k(t)$ in $L_2([0,T])$ given by
\begin{equation}
h_1(t) = \frac{1}{\sqrt{T}},~~h_k(t) = \sqrt{\frac{2}{T}}\cos\left(\frac{(k-1)\pi t}{T}\right), ~~k \in \{2,3,\ldots\}.
\end{equation}
Let us first study the effect of local improvements to POD basis functions. Choose $\bm{\hat{\xi}}$ and $\bm{\hat{\zeta}}$ to be two $d$ dimensional vectors with Gaussian distributed randomly chosen entries. We again consider the transformation $\bm{\hat{\xi}} + \theta(\bm{\hat{\zeta}}-\bm{\hat{\xi}})/||\bm{\hat{\zeta}}-\bm{\hat{\xi}}||$. We let $\theta$ to change in the interval $[-1/2,1/2]$. We compute the POD basis functions at $\bm{\hat{\xi}}$ and use them, along with their extrapolated and expanded local improvements, to solve the nearby problem at $\bm{\hat{\xi}} + \theta(\bm{\hat{\zeta}}-\bm{\hat{\xi}})/||\bm{\hat{\zeta}}-\bm{\hat{\xi}}||$ using the Group POD method. We then compare the results with the ones from the GFE method. For the finite element code we partition the spatial domain $[0,1]$ into 64 intervals. We also divide the time domain $[0,T]$ into 200 time steps. For the reduced order models, we use 10 POD basis functions which results in 10 extrapolated basis functions and 10+10 extended basis functions. We also consider the case when 20 POD basis functions and consequently 20 extrapolated bases are employed, for comparison. We let $d=10$.
\begin{figure}[h!]
\centering\includegraphics[width=100mm]{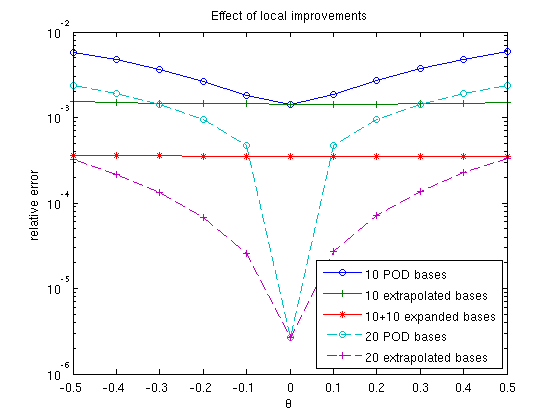}
\caption{Comparison of the effects of extrapolated and expanded local improvements to POD basis functions}
\label{Sensitivity}
\end{figure}
Figure \ref{Sensitivity} shows the typical behaviour of Group POD method when POD, Extrapolated, and Expanded basis functions are employed. The authors believe that the improvement in accuracy that is achieved by employing Expanded basis functions is not worth the increase in the number of bases. Therefore, in the multi-fidelity collocation method proposed earlier we use Extrapolated basis functions to locally improve the performance of Group POD method.

\begin{figure}[h!]
\centering\includegraphics[width=100mm]{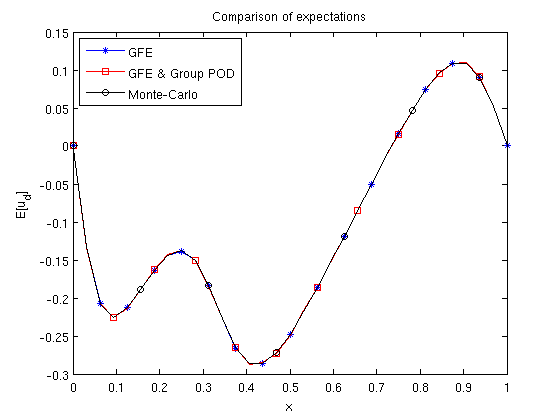}
\centering\includegraphics[width=100mm]{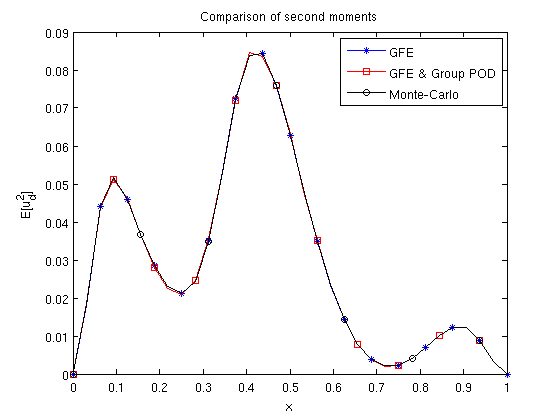}
\caption{Expectation and Second Moment of the solution given at the final time $T=0.8$ when $d=3$, for sparse grid level 8 and $\eta=0.5$.}
\label{Comparison}
\end{figure}

Now we apply the multi-fidelity collocation method. Let us assume that $\eta = 1/2$ and $d=3$. For the finite element code we partition the spatial domain $[0,1]$ into 32 intervals. We also divide the time domain $[0,T]$ into 20 time steps. We use $10$ extrapolated POD basis functions. We employ the Smolyak algorithm with sparse grid level $q=8$. We use the Clenshaw-Curtis abscissas (see \cite{clenshaw1960method}) as collocation points. These abscissas are the extrema of Chebyshev polynomials. In figure \ref{Comparison} we compare expectations and second moments of solutions resulting from our multi-fidelity method (GFE \& Group POD) and the ones resulting from a regular sparse grid method with the full employment of GFE as the high fidelity algorithm. In this figure, we are also including the solutions coming out of the Monte-Carlo method with the full employment of the high fidelity algorithm for reference.
For this value of $\eta = 1/2$, the number of times that our algorithm calls the high fidelity (GFE) code is reduced to $2030$. Compare it to $6018$, the number of times that the GFE code is called when sparse grid stochastic collocation with full employment of the high-fidelity algorithm is utilized. Comparing the two methods, assuming that sparse grid stochastic collocation with full employment of the high-fidelity algorithm is accurate enough, we get that the error in expectation is given by $8.2\times 10^{-4}$ and the one in the second moment is $1.4\times 10^{-3}$.

\begin{figure}[h!]
\centering\includegraphics[width=100mm]{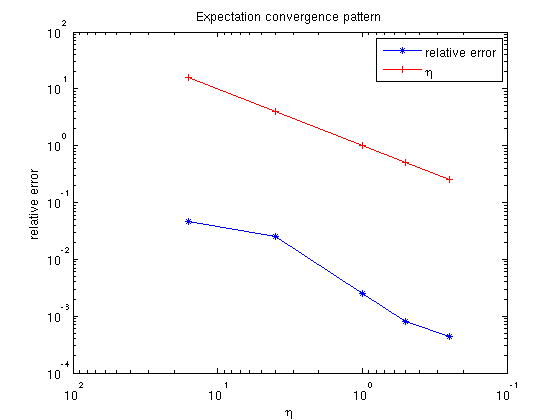}
\centering\includegraphics[width=100mm]{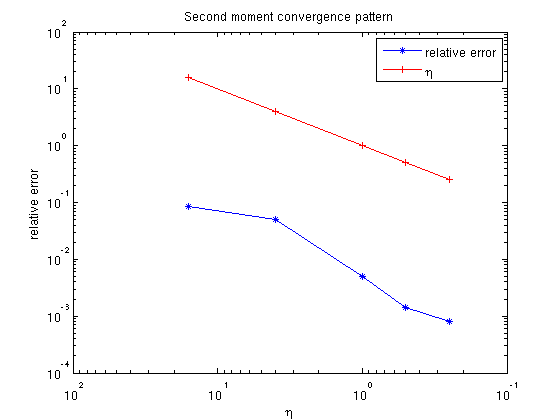}
\caption{$L_2([0,1])$-norm convergence pattern of Expectation and Second Moment of the solution given at the final time $T=0.8$ when $d=3$, for sparse grid level 8.}
\label{ConvergencePatern}
\end{figure}

Figure \ref{ConvergencePatern}, shows the convergence patterns of expectations and second moments of solutions with regard to $\eta$. We are in fact comparing our multi-fidelity method with a regular sparse grid stochastic collocation method. Note that for small enough $\eta$ (less than the shortest distance between the collocation points) we get the regular sparse grid method back. Therefore the error is zero for such a small $\eta$.

\begin{figure}[h!]
\centering\includegraphics[width=100mm]{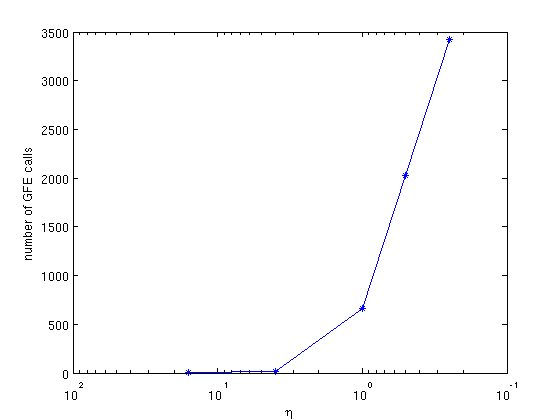}
\caption{Number of times that the high fidelity algorithm (GFE method) is called as a function of $\eta$.}
\label{GFECalls}
\end{figure}

Figure \ref{GFECalls} demonstrates how the number of times that the finite element code is employed increases with respect to a decrease in $\eta$.

Table \ref{Table} is just another way of presenting the data depicted in figures \ref{ConvergencePatern} and \ref{GFECalls}.
\begin{table}[h!]
\caption{Relative errors and the number of times that the finite element code is employed for different values of $\eta$.}
\begin{center}
\label{Table}
\begin{tabular}{|c|c|c|c|c|c|}
\hline 
$\eta$ & $\#$ FE calls & Expectation $L_2$ error & Variance $L_2$ error \\ 
\hline 
16 &	1 &	4.56E-02 &	8.60E-02\\
\hline
4 &	20 &	2.53E-02 &	4.94E-02\\
\hline
1 &	659 &	2.55E-03 &	4.95E-03\\
\hline
$1/2$ &	 2030 &	8.16E-04 &	1.43E-03\\
\hline
$(1/4)$ &	3424	 & 4.45E-04 &	8.05E-04\\
\hline
\end{tabular} 
\end{center}
\end{table}

\begin{remark} 
The method proposed in this work is a generalization of the one introduced in \cite{maziar1}. It also generalizes the one in \cite{maziar2}, where rigorous convergence analysis of this method is provided for parabolic partial differential equations driven by random input data.
\end{remark}
\section{Conclusion and discussion}
The main aim of this work is to show how sparse grid stochastic collocation method can be enhanced by model reduction techniques and by sensitivity analysis of Proper Orthogonal Decomposition basis functions. We have chosen to apply it to the stochastic Burgers equation because of its tractability and because lots of powerful model reduction techniques such as group POD existed for this equation. Our method enhances the power of stochastic collocation methods in handling more stochastic differential equations. Our paper provides another reason for the importance of research in Model reduction techniques. We believe that more research can follow our work and apply our method to more interesting problems. Such problems can be in fluid dynamics and fluid structure interactions. Navier-Stokes equations under uncertainty can for instance be a reasonable next step.
\section*{Acknowledgments}
The authors thank George Karniadakis, Dongbin Xiu, Alireza Doostan, Sergey Lototsky and Peter Kloeden for their valuable comments during the ICERM uncertainty quantification workshop at Brown university. We also thank Karen Willcox for her comments during her visit to George Mason University.

\bibliographystyle{siam}
\bibliography{mybib}{}

\end{document}